\documentclass[dvipdfmx]{article}
\usepackage{amsmath,amssymb,amsfonts,amsthm}
\usepackage{mathtools}
\usepackage{mathrsfs}
\usepackage{tasks}
\usepackage[inline]{enumitem}
\usepackage[hidelinks]{hyperref} 
\usepackage[dvipdfmx]{graphicx,xcolor}
\usepackage{tikz}
\usepackage{xparse} 
\usepackage{xspace} 
\usepackage{setspace}
\usepackage{bussproofs}

\theoremstyle{plain}
\newtheorem{theorem}{Theorem}[section]
\newtheorem*{theorem*}{Theorem}
\newtheorem{proposition}[theorem]{Proposition}
\newtheorem*{proposition*}{Proposition}

\newtheorem*{lemma*}{Lemma}

\newtheorem*{claim*}{Claim}
\newtheorem{corollary}[theorem]{Corollary}
\newtheorem*{corollary*}{Corollary}
\newtheorem{fact}[theorem]{Fact}
\newtheorem*{fact*}{Fact}

\newtheorem*{guide*}{Guide}

\theoremstyle{definition}
\newtheorem{definition}[theorem]{Definition}

\newtheorem{example}[theorem]{Example}
\newtheorem*{example*}{Example}
\newtheorem{remark}[theorem]{Remark}
\newtheorem*{remark*}{Remark}

\DeclarePairedDelimiterX{\set}[1]{\{}{\}}{\setargs{#1}}
\NewDocumentCommand{\setargs}{>{\SplitArgument{1}{;}}m}
{\setargsaux#1}
\NewDocumentCommand{\setargsaux}{mm}
{\IfNoValueTF{#2}{#1} {#1\nonscript\:\delimsize\vert\allowbreak\nonscript\:\mathopen{}#2}}%
\def\Set{\set*}

\newcommand*{\tif}{\mathrel{\text{if}}}
\newcommand*{\tnot}{\mathrel{\text{not}}}
\newcommand*{\tand}{\mathrel{\text{and}}}

\newcommand*{\PropVar}{\mathrm{PropVar}}
\newcommand*{\CL}{\mathbf{Cl}}
\newcommand*{\Lang}{\mathscr{L}}
\newcommand*{\PF}{\Lang_\mathrm{p}}

\newcommand*{\K}{\mathbf{K}}
\newcommand*{\N}{\mathbf{N}}
\newcommand*{\NA}[2]{\mathbf{NA}_{#1,#2}}
\newcommand*{\NpA}[2]{\mathbf{N}^+\mathbf{A}_{#1,#2}}
\newcommand*{\NRA}[2]{\mathbf{NRA}_{#1,#2}}
\newcommand*{\NxA}[2]{\mathbf{N}^\ast\mathbf{A}_{#1,#2}}
\newcommand*{\NAmn}{\NA{m}{n}}
\newcommand*{\NpAmn}{\NpA{m}{n}}
\newcommand*{\NRAmn}{\NRA{m}{n}}
\newcommand*{\NxAmn}{\NxA{m}{n}}
\newcommand*{\Amn}[1][m,n]{\mathrm{Acc}_{#1}}
\newcommand*{\Ros}{\textsc{Ros}}
\newcommand*{\RosBox}{\text{\textsc{Ros}}^{\Box}}
\newcommand*{\FMT}{Fitting, Marek, and Truszczy\'{n}ski\xspace}
\newcommand*{\MF}{\Lang_\Box}

\newcommand*{\LK}{\mathbf{LK}}
\newcommand*{\GNA}[2]{\mathbf{GNA}_{#1,#2}}
\newcommand*{\GNpA}[2]{\mathbf{GN}^+\mathbf{A}_{#1,#2}}
\newcommand*{\GNRA}[2]{\mathbf{GNRA}_{#1,#2}}
\newcommand*{\GNxA}[2]{\mathbf{GN}^\ast\mathbf{A}_{#1,#2}}
\newcommand*{\GNAmn}{\GNA{m}{n}}
\newcommand*{\GNpAmn}{\GNpA{m}{n}}
\newcommand*{\GNRAmn}{\GNRA{m}{n}}
\newcommand*{\GNxAmn}{\GNxA{m}{n}}
\newcommand*{\To}{\mathrel{\Longrightarrow}}
\newcommand*{\SequentRuleName}[2]{\ensuremath{\mathrm{\mathord{#1}#2}}}
\newcommand*{\Var}{\mathrm{V}}
\newcommand*{\VarPos}{\mathrm{V}^+}
\newcommand*{\VarNeg}{\mathrm{V}^-}
\newcommand*{\VarAny}{\mathrm{V}^\bullet}
\newcommand*{\VarAnother}{\mathrm{V}^\circ}

\usepackage{subfiles}

\title{Uniform Lyndon interpolation for the pure logic of necessitation with a modal reduction principle}
\author{Yuta Sato \footnote{Email: 231x032x@gsuite.kobe-u.ac.jp} \footnote{Graduate School of System Informatics,
Kobe University,
1-1 Rokkodai, Nada, Kobe 657-8501, Japan.}}
\date{}

\begin{document}

\maketitle

\begin{abstract}
We prove the uniform Lyndon interpolation property (ULIP) of some extensions of the pure logic of necessitation $\N$.
For any $m, n \in \mathbb{N}$, $\NpAmn$ is the logic obtained from $\N$ by adding a single axiom $\Box^n \varphi \to \Box^m \varphi$, $\Diamond$-free modal reduction principle,
together with a rule $\frac{\neg \Box \varphi}{\neg \Box \Box \varphi}$, required to make the logic complete with respect to its Kripke-like semantics.
We first introduce a sequent calculus $\GNpAmn$ for $\NpAmn$ and show that it enjoys cut elimination, proving Craig and Lyndon interpolation properties as consequences.
We then introduce a general method, called propositionalization, that allows one to prove ULIP of a logic by reducing it to that of some weaker logic.
Lastly, we construct a propositionalization of $\NpAmn$ into classical propositional logic $\CL$, proving ULIP as a corollary.
We also prove ULIP of $\NAmn = \N + \Box^n \varphi \to \Box^m \varphi$ and $\NRAmn = \N + \Box^n \varphi \to \Box^m \varphi + \frac{\neg \varphi}{\neg \Box \varphi}$ in the same manner.
\end{abstract}

\def\fCenter{\ \To\ }

\section{Introduction}

The pure logic of necessitation $\N$, originally introduced by \FMT \cite{FMT92},
is a nonnormal modal logic obtained from classical propositional logic $\CL$ by adding the necessitation rule $\frac{\varphi}{\Box \varphi}$,
or from the modal logic $\K$ by removing the $\mathrm{K}$ axiom.
They introduced Kripke-like relational semantics for $\N$, namely $\N$-frames and $\N$-models,
and proved the finite frame property (FFP) of $\N$ with respect to the class of all $\N$-frames\footnote{
  It is worth noting that Omori and Skurt \cite{OS24} recently rediscovered the same logic as $\N$, namely $\mathbf{M}^+$ in their paper,
  and gave a nondeterministic, many-valued semantics.
}.
They then used FFP of $\N$ to analyze nonmonotonic reasoning, and also suggested that the $\Box$-modality
in $\N$ can be interpreted with their notion of provability, reading $\Box \varphi$ as \emph{$\varphi$ is already derived}.

Decades later, Kurahashi \cite{Kur23} discovered that $\N$ and its extensions play an important role in the context of provability logic,
which uses more traditional notion of provability than Fitting et al.\ employed.
Kurahashi proved FFP of several extensions of $\N$ including $\mathbf{N4} = \N + \Box \varphi \to \Box \Box \varphi$,
and showed that $\N$ is exactly the provability logic of all provability predicates, and that $\mathbf{N4}$ is exactly the provability logic of all provability predicates satisfying an additional condition called \textsf{D3}.
Generalizing Kurahashi's work on FFP of $\mathbf{N4}$, Kurahashi and the author \cite{KS24} proved that:
\begin{itemize}
  \item The logic $\NAmn = \N + \Box^n \varphi \to \Box^m \varphi$ is not complete for some combinations of $m, n \in \mathbb{N}$, and;
  \item The logic $\NpAmn = \NAmn + \frac{\neg \Box \varphi}{\neg \Box \Box \varphi}$ has FFP for any $m, n \in \mathbb{N}$.
\end{itemize}
We shall refer to $\Box^n \varphi \to \Box^m \varphi$ as the $\Amn$ axiom.
These results on completeness and FFP are interesting when compared to those of their normal counterpart $\K + \Amn$.
There is a famous result by Lemmon and Scott\,\cite{LS77} that $\K + \Amn$ is complete for every $m, n \in \mathbb{N}$,
while it is still unknown to this day whether $\K + \Amn$ has FFP in general (cf.~\cite[Problem 12.1]{CZ97}, \cite[Problem 6]{WZ07}, \cite{Zak97} and \cite[Problem 6.12]{Zak97_2}). 
These results highlight intriguing differences between the extensions of $\N$ and their normal counterparts.

In this paper, we study Craig (CIP), Lyndon (LIP), and uniform (UIP) interpolation properties of $\NpAmn$,
which was proposed in \cite[Problem 3]{KK23} and \cite{KS24} as open problems,
to further investigate the difference between $\NpAmn$ and $\K + \Amn$\footnote{For a comprehensive survey of interpolation results in superintuitionistic and modal logics, we refer to Gabbay and Maksimova \cite{GM09}.}.
Not much is known about the interpolation properties of $\K + \Amn$ in general,
but several interesting partial results were obtained.
Gabbay \cite{Gab72} proved that $\K + \Amn[n, 1]$ enjoys CIP for $n > 0$ and
Kuznets \cite{Kuz16} extended Gabbay's result to prove LIP of it.
For negative results, Bílková \cite{Bil07} proved that $\mathbf{K4}$ lacks UIP based on the failure of UIP in $\mathbf{S4}$ due to Ghilardi and Zawadowski \cite{GZ95},
while Marx \cite{Mar95} proved that $\K + \Amn[1, 2]$ lacks even CIP, yet alone LIP and UIP.
Given these results, it would be natural to ask what happens to the interpolation properties of $\mathbf{K4}$ and $\K + \Amn[1, 2]$ if the $\mathrm{K}$ axiom is removed,
and we shall show that, as a consequence on the general result on $\NpAmn$, both $\mathbf{N4}$ and $\N + \Amn[1, 2]$ enjoy all of CIP, LIP, and UIP,
which suggests that the problem lies not in the $\Amn$ axiom, but rather in the combination of $\Amn$ and $\mathrm{K}$.

In addition to $\NpAmn$, we also consider two more logics in this paper,
$\NAmn$, which is incomplete for some, but still sound for all, combinations of $m, n \in \mathbb{N}$,
and $\NRAmn = \NAmn + \frac{\neg \varphi}{\neg \Box \varphi}$,
which is a generalization of the logic $\mathbf{NR4}$ introduced by Kurahashi \cite{Kur23}.
For the sake of brevity, we will treat all three logics at once in the form $\NxAmn \in \set{ \NAmn,\, \NpAmn,\, \NRAmn }$.

This paper is organized as follows. We first present some preliminary definitions and results in Section \ref{sec:prelim}.
The uniform Lyndon interpolation property (ULIP), a strengthening of both UIP and LIP, introduced by Kurahashi \cite{Kur20},
is also presented in this section along with some facts on the relationship with UIP and LIP.
Then in Section \ref{sec:sequent}, we introduce a sequent calculus $\GNxAmn$ for $\NxAmn$ and prove that it enjoys cut elimination.
As a consequence of the cut elimination theorem, we first prove CIP and LIP of $\NxAmn$ with Maehara's method.
Then in Section \ref{sec:p18n}, we introduce a general method called \emph{propositionalization}, stating that
ULIP of a weaker logic $L$ can be transferred to that of a stronger logic $M$ if there is an embedding from $M$ to $L$ with certain properties.
Lastly in Section \ref{sec:ulip}, we prove that $\NxAmn$ enjoys ULIP by constructing such an embedding.
We first define a pair of translations from modal formulae to propositional formulae, carefully constructed to include the information needed to emulate the rules of $\GNxAmn$ in $\LK$.
We then prove that a cut-free proof in $\GNxAmn$ can indeed be emulated in $\LK$ using the said translations,
establishing an embedding from $\NxAmn$ to $\CL$.
We then show that the embedding satisfies the required properties of propositionalization.
As a consequence, we obtain ULIP of $\NxAmn$.

\section{Preliminaries} \label{sec:prelim}

Let $\PropVar$ denote the set of all propositional variables.
Let $\PF = \set{ \land, \lor, \to, \bot }$ be the language of $\CL$ and $\MF = \set{ \land, \lor, \to, \bot, \Box }$ be the language of modal logics.
We use $\top$ and $\neg \varphi$ as abbreviations of $\bot \to \bot$ and $\varphi \to \bot$, respectively.

\subsection{The logics $\NAmn$, $\NpAmn$, and $\NRAmn$}

Let $m, n \in \mathbb{N}$.
The logic $\N$ is obtained from classical propositional logic $\CL$ by adding the rule $\textsc{Nec}$: $\frac{\varphi}{\Box \varphi}$ with extending the language $\Lang_p$ to $\Lang_\Box$,
$\NAmn$ is obtained from $\N$ by adding the axiom $\Amn$: $\Box^n \varphi \to \Box^m \varphi$,
$\NpAmn$ is obtained from $\NAmn$ by adding the rule $\RosBox$: $\frac{\neg \Box \varphi}{\neg \Box \Box \varphi}$,
and $\NRAmn$ is obtained from $\NAmn$ by adding an even stronger rule $\Ros$: $\frac{\neg \varphi}{\neg \Box \varphi}$
\footnote{$\Amn$ stands for $(m,n)$-accessibility. $\Ros$ originates from the derivability condition of Rosser provability predicates (cf.\ \cite{Kur23}). $\RosBox$ is a weakening of $\Ros$ with an additional box.}.
FFP of $\NAmn$ and $\NpAmn$ are studied by Kurahashi and the author \cite{KS24}.

\begin{proposition}[{\cite[Section 4]{KS24}}] \leavevmode
  \begin{itemize}
    \item $\RosBox$ is admissible in every normal modal logic, including $\K + \Amn$.
    \item $\RosBox$ is also admissible in $\NAmn$ (i.e.\ $\NpAmn = \NAmn$) when $m > 0$ or $n < 2$.
    \item $\NpAmn \supsetneq \NAmn$ when $m = 0$ and $n \ge 2$.
  \end{itemize}
  Consequently, $\K + \Amn \supsetneq \NpAmn \supseteq \NAmn$.
\end{proposition}

\emph{$\N$-frames}, the Kripke-like semantics for $\N$, are introduced in \cite{FMT92}.
Several frame conditions on $\N$-frames, namely transitivity, seriality, and \emph{$(m,n)$-accessibility}, are presented in \cite{Kur23} and \cite{KS24}.
Here, $(m,n)$-accessibility (``if one can see the other in $m$ steps, then can also see in $n$ steps'') is a natural generalization of transitivity,
which corresponds to the fact that $\Amn$: $\Box^n \varphi \to \Box^m \varphi$ generalizes $\mathrm{4}$: $\Box \varphi \to \Box \Box \varphi$.

\begin{theorem}[{\cite[Corollary 4.6, Theorem 5.1]{KS24}}] \label{thm:nramn-ffp}
  With respect to the class of all $(m,n)$-accessible $\N$-frames,
  \begin{itemize}
    \item $\NAmn$ is incomplete when $m = 0$ and $n \ge 2$.
    \item $\NpAmn$ enjoys FFP for every $m, n \in \mathbb{N}$.
  \end{itemize}
\end{theorem}

These results show that, although it is admissible in many cases, $\RosBox$ is an important rule for extensions of $\N$ in general,
with respect to the relational semantics by Fitting et al.

$\NRAmn$ is a generalization of the logic $\mathbf{NR4} = \mathbf{N4} + \Ros$, which was introduced and analyzed by Kurahashi \cite{Kur23}.
The completeness and FFP of $\NRAmn$, to the best of the author's knowledge, are not yet known.
The soundness of it, however, can easily be obtained.

\begin{proposition} \label{prop:nramn-sound}
  $\NRAmn$ is sound with respect to the class of all $(m,n)$-accessible and serial $\N$-frames.
  \begin{proof}
    Follows from \cite[Theorem 3.3]{KS24} and \cite[Corollary 3.7]{Kur23}.
  \end{proof}
\end{proposition}

\subsection{A handful of interpolation properties}

\begin{definition}
  Let $\mathscr{L} \in \set{ \PF, \MF }$.
  We inductively define the sets of positive and negative variables of $\varphi \in \mathscr{L}$,
  denoted by $\VarPos(\varphi)$ and $\VarNeg(\varphi)$ respectively:
  \begin{itemize}
    \item $\VarPos(p) = \set{p}$ and $\VarNeg(p) = \VarPos(\bot) = \VarNeg(\bot) = \emptyset$;
    \item $\VarAny(\psi_1 \to \psi_2) = \VarAnother(\psi_1) \cup \VarAny(\psi_2)$, for each $(\bullet, \circ) \in \set{(+, -), (-, +)}$;
    \item $\VarAny(\psi_1 \circledcirc \psi_2) = \VarAny(\psi_1) \cup \VarAny(\psi_2)$, for each $\bullet \in \set{+, -}$ and $\circledcirc \in \set{\land, \lor}$;
    \item $\VarAny(\Box \psi) = \VarAny(\psi)$, for each $\bullet \in \set{+, -}$.
  \end{itemize}
  And for any set $\Gamma \subseteq \mathscr{L}$, we let $\VarAny(\Gamma) = \bigcup_{\varphi \in \Gamma} \VarAny(\varphi)$, for each $\bullet \in \set{+, -}$.
  We also let $\Var(\varphi) = \VarPos(\varphi) \cup \VarNeg(\varphi)$ and $\Var(\Gamma) = \VarPos(\Gamma) \cup \VarNeg(\Gamma)$.
\end{definition}

\begin{definition} \label{def:cip-lip}
  A logic $L$ with a language $\mathscr{L}$ is said to enjoy \emph{Lyndon interpolation property (LIP)}
  if for every $\varphi, \psi \in \mathscr{L}$ such that $L \vdash \varphi \to \psi$,
  there is $\chi \in \mathscr{L}$ that satisfies the following conditions:
  \begin{enumerate}
    \item $\VarAny(\chi) \subseteq \VarAny(\varphi) \cap \VarAny(\psi)$, for each $\bullet \in \set{+, -}$;
    \item $L \vdash \varphi \to \chi$ and $L \vdash \chi \to \psi$.
  \end{enumerate}
  Such $\chi$ is called an \emph{interpolant} of $\varphi \to \psi$.
  A logic $L$ is said to enjoy \emph{Craig interpolation property (CIP)} if it satisfies the above,
  but with the first condition of $\chi$ replaced by $\Var(\chi) \subseteq \Var(\varphi) \cap \Var(\psi)$.
\end{definition}

\begin{definition} \label{def:ulip-uip}
  A logic $L$ with a language $\mathscr{L}$ is said to enjoy \emph{Uniform Lyndon interpolation property (ULIP)} if
  for any $\varphi \in \mathscr{L}$ and any finite $P^\bullet \subseteq \PropVar$ for each $\bullet \in \set{+, -}$,
  there is $\chi \in \mathscr{L}$ such that:
  \begin{enumerate}
    \item $\VarAny(\chi) \subseteq \VarAny(\varphi) \setminus P^\bullet$, for each $\bullet \in \set{+, -}$;
    \item $L \vdash \varphi \to \chi$;
    \item $L \vdash \chi \to \psi$ for any $\psi \in \mathscr{L}$ such that $L \vdash \varphi \to \psi$ and $\VarAny(\psi) \cap P^\bullet = \emptyset$, for each $\bullet \in \set{+, -}$.
  \end{enumerate}
  Such $\chi$ is called a \emph{post-interpolant} of $(\varphi, P^+, P^-)$.
  A logic $L$ is said to enjoy \emph{Uniform interpolation property (UIP)} if it satisfies the above,
  but with all occurrences of the superscript $\bullet$ and the condition ``for each $\bullet \in \set{+, -}$'' omitted everywhere.
\end{definition}

Here, CIP easily follows from LIP.
It is also easy to see that UIP implies CIP by taking any $\psi$ such that $L \vdash \varphi \to \psi$ and letting $P = \Var(\varphi) \setminus \Var(\psi)$.
The fact that ULIP implies UIP and LIP is similarly proven in \cite{Kur20}.

\subsection{Sequent calculi and Maehara's method}

CIP and LIP of a logic can be proven in many ways, syntactically, semantically, or algebraically.
One famous way to do this is \emph{Maehara's method}, which utilizes a sequent calculus for the logic that enjoys cut elimination.

\begin{definition}
  Let $\mathscr{L}$ be any language.
  Let $\varphi, \psi, \ldots \in \mathscr{L}$, and let $\Gamma, \Delta, \ldots$ range over finite subsets of $\mathscr{L}$.
  \begin{itemize}
    \item A pair $(\Gamma, \Delta)$ is called a \emph{sequent} and is written as $\Gamma \To \Delta$.
    \item We write $\varphi, \Gamma$ to mean $\set{\varphi} \cup \Gamma$, and $\Gamma, \Delta$ to mean $\Gamma \cup \Delta$.
          We also permit the omission of emptysets.
          With this in mind, we write
          $\set{\varphi_1, \varphi_2} \To \set{\psi_1, \psi_2}$ as $\varphi_1, \varphi_2 \To \psi_1, \psi_2$,
          write $\emptyset \To \set{\varphi}$ as $\To \varphi$, and so on.
    \item We write $(\Gamma \To \Delta) = (\Gamma' \To \Delta')$ to mean $\Gamma = \Gamma'$ and $\Delta = \Delta'$.
          When we write $(\Gamma \To \Delta) = (\varphi, \Gamma' \To \Delta', \psi)$, we assume that
          $\varphi \notin \Gamma'$ and $\psi \notin \Delta'$.
    \item A pair of sequents $(\Gamma_1 \To \Delta_1; \Gamma_2 \To \Delta_2)$ is called a \emph{partition} of a sequent $\Gamma \To \Delta$
          if $\Gamma_1 \cup \Gamma_2 = \Gamma$, $\Delta_1 \cup \Delta_2 = \Delta$,
          and $\Gamma_1 \cap \Gamma_2 = \Delta_1 \cap \Delta_2 = \emptyset$.
  \end{itemize}
\end{definition}

\begin{definition} \label{def:LK}
  The propositional calculus $\LK$ has the following initial sequents and rules.

  \subparagraph*{Initial sequents}
  Initial sequents have the form of $\varphi \To \varphi$ or $\bot \To$.

  \subparagraph*{Logical rules}
  \begin{center}
    \Axiom$\Gamma, \varphi_i \fCenter \Delta$
    \RightLabel{(\SequentRuleName{\land}{L}, $i = 1,2$)}
    \UnaryInf$\Gamma, \varphi_1 \land \varphi_2 \fCenter \Delta$
    \DisplayProof
    \quad
    \Axiom$\Gamma \fCenter \varphi_1, \Delta$
    \Axiom$\Gamma \fCenter \varphi_2, \Delta$
    \RightLabel{(\SequentRuleName{\land}{R})}
    \BinaryInf$\Gamma \fCenter \varphi_1 \land \varphi_2, \Delta$
    \DisplayProof
  \end{center}

  \begin{center}
    \Axiom$\Gamma, \varphi_1 \fCenter \Delta$
    \Axiom$\Gamma, \varphi_2 \fCenter \Delta$
    \RightLabel{(\SequentRuleName{\lor}{L})}
    \BinaryInf$\Gamma, \varphi_1 \lor \varphi_2 \fCenter \Delta$
    \DisplayProof
    \quad
    \Axiom$\Gamma \fCenter \varphi_i, \Delta$
    \RightLabel{(\SequentRuleName{\lor}{R}, $i = 1,2$)}
    \UnaryInf$\Gamma \fCenter \varphi_1 \lor \varphi_2, \Delta$
    \DisplayProof
  \end{center}

  \begin{center}
    \Axiom$\Gamma_1 \fCenter \varphi, \Delta_1$
    \Axiom$\Gamma_2, \psi \fCenter \Delta_2$
    \RightLabel{(\SequentRuleName{\to}{L})}
    \BinaryInf$\Gamma_1, \Gamma_2, \varphi \to \psi \fCenter \Delta_1, \Delta_2$
    \DisplayProof
    \quad
    \Axiom$\Gamma, \varphi \fCenter \psi, \Delta$
    \RightLabel{(\SequentRuleName{\to}{R})}
    \UnaryInf$\Gamma \fCenter \varphi \to \psi, \Delta$
    \DisplayProof
  \end{center}

  \subparagraph*{Structural rules}

  \begin{center}
    \Axiom$\Gamma \fCenter \Delta$
    \RightLabel{(wL)}
    \UnaryInf$\Gamma, \varphi \fCenter \Delta$
    \DisplayProof
    \quad
    \Axiom$\Gamma \fCenter \Delta$
    \RightLabel{(wR)}
    \UnaryInf$\Gamma \fCenter \varphi, \Delta$
    \DisplayProof
  \end{center}
  In all logical and structural rules, we call the newly obtained formula (e.g. $\varphi_1 \land \varphi_2$ in \SequentRuleName{\land}{L})
  the \emph{principal formula} of the rule.

  \subparagraph*{Cut rule}
  \begin{center}
    \Axiom$\Gamma_1 \fCenter \Delta_1, \varphi$
    \Axiom$\varphi, \Gamma_2 \fCenter \Delta_2$
    \RightLabel{(cut)}
    \BinaryInf$\Gamma_1, \Gamma_2 \fCenter \Delta_1, \Delta_2$
    \DisplayProof
  \end{center}
  The formula $\varphi$ that is `cut' in the cut rule is called the \emph{cut formula}.
\end{definition}

\begin{definition}
  Let $L$ be a logic and $C$ be a sequent calculus.
  We say $C$ is \emph{a sequent calculus for $L$} (or $L$ is \emph{with a sequent calculus $C$}) if
  $(C \vdash \Gamma \To \Delta)$ is equivalent to $\left(L \vdash \bigwedge \Gamma \to \bigvee \Delta\right)$.
\end{definition}

\begin{fact}[cf. \cite{Ono19}] \label{prop:lk-pc}
  $\LK$ is a sequent calculus for $\CL$.
\end{fact}

\begin{definition}
  A sequent calculus $C$ is said to enjoy \emph{cut elimination}
  if for every sequent $\Gamma \To \Delta$ that is provable in $C$,
  there is a proof of $\Gamma \To \Delta$ in $C$
  in which the cut rule is not used. Such a proof is called a \emph{cut-free proof}.
\end{definition}

\begin{fact}[cf. \cite{Ono19}] \label{thm:lk-cut-elim}
  $\LK$ enjoys cut elimination.
\end{fact}

Now for a logic $L$ with a sequent calculus $C$ that enjoys cut elimination,
one can apply Maehara's method as follows:
Take any sequent $\Gamma \To \Delta$ that is provable in $C$,
then we use induction on the length of the cut-free proof to
show that for any partition $(\Gamma_1 \To \Delta_1; \Gamma_2 \To \Delta_2)$ of $\Gamma \To \Delta$,
there is $\chi$ such that
both $\Gamma_1 \To \Delta_1, \chi$ and $\chi, \Gamma_2 \To \Delta_2$ are provable in $C$, and
$\VarAny(\chi) \subseteq (\VarAny(\Gamma_1) \cup \VarAnother(\Delta_1)) \cap (\VarAnother(\Gamma_2) \cup \VarAny(\Delta_2))$ for each $(\bullet, \circ) \in \set{(+, -), (-, +)}$.
Then, LIP of $L$ is a direct consequence of it.

\begin{fact}[\cite{Lyn59}]
  $\CL$ enjoys LIP and, consequently, also enjoys CIP.
\end{fact}

\subsection{Several ways to prove ULIP}

U(L)IP is easily proven for a logic with a property called \emph{local tabularity}:

\begin{definition}
  Let $L$ be a logic with a language $\mathscr{L}$.
  \begin{itemize}
    \item For $\varphi, \psi \in \mathscr{L}$, let us write $\varphi \equiv_L \psi$ to mean $L \vdash \varphi \leftrightarrow \psi$, then $\equiv_L$ is clearly an equivalence relation.
    \item $L$ is said to be \emph{locally tabular} if for any finite $P \subseteq \PropVar$,
          the set $T(P) = \set*{ \varphi; \Var(\varphi) \subseteq P } / \mathord{\equiv_L}$ is finite.
  \end{itemize}
\end{definition}

\begin{proposition}[\cite{Mak14}, \cite{Kur20}] \leavevmode
  For a locally tabular logic $L$ with a language $\mathscr{L}$,
  \begin{enumerate}
    \item[(1)] If $L$ enjoys CIP, then it also enjoys UIP.
    \item[(2)] If $L$ enjoys LIP, then it also enjoys ULIP.
  \end{enumerate}
  \begin{proof}
    We shall only prove (1), as (2) is similarly proven.
    Take any $\varphi \in \mathscr{L}$ and any finite $P \subseteq \PropVar$.
    Pick one representative $\chi_i$ of $[\chi_i]$
    for each $[\chi_i] \in T(\Var(\varphi) \setminus P)$ such that $L \vdash \varphi \to \chi_i$.
    Let $X$ be the set of every such $\chi_i$, then $X$ is finite by the local tabularity of $L$.
    We let $\hat{\chi} = \bigwedge X$, and take any $\psi \in \mathscr{L}$ such that $L \vdash \varphi \to \psi$ and $\Var(\psi) \cap P = \emptyset$.
    Then by CIP of $L$, there is $\chi \in \mathscr{L}$ such that $L \vdash \varphi \to \chi$, $L \vdash \chi \to \psi$, and $\Var(\chi) \subseteq \Var(\varphi) \cap \Var(\psi) \subseteq \Var(\varphi) \setminus P$.
    Here, there is $\chi_i \in X$ such that $L \vdash \chi_i \leftrightarrow \chi$, then $L \vdash \hat{\chi} \to \chi$, so $L \vdash \hat{\chi} \to \psi$.
    Therefore, $\hat{\chi}$ is the post-interpolant of $(\varphi, P^+, P^-)$.
  \end{proof}
\end{proposition}

\begin{fact}[cf. \cite{CZ97}]
  $\CL$ is locally tabular and, consequently, enjoys ULIP.
\end{fact}

\begin{fact}
  For $L \in \set{ \N,\, \NAmn,\, \NpAmn,\, \NRAmn }$, $L$ is not locally tabular.
  \begin{proof}
    It is easily proven that $\Box \neg^{2k} p \not \equiv_L \Box \neg^{2l} p$ for any $k \ne l$
    using the relational semantics by \cite{FMT92}, \cite{Kur23}, and \cite{KS24}.
  \end{proof}
\end{fact}

Several methods for proving U(L)IP of logics without local tabularity are also proposed.
Visser \cite{Vis96} proved UIP of some normal modal logics semantically by bounded bisimulation,
and Kurahashi \cite{Kur20} extended Visser's method to prove ULIP of several modal logics.
Iemhoff \cite{Iem19} proved that a logic enjoys UIP if it has a sequent calculus with certain structural properties,
and then gave a syntactic proof of UIP of several intuitionistic modal logics.
Akbar Tabatabai et al.\ \cite{Tab24} employed a similar method to prove ULIP of several nonnormal modal logics and conditional logics.

A carefully constructed translation enables one to prove U(L)IP of a logic by reducing it to that of another logic.
Visser \cite{Vis96} proved that UIP of intuitionistic propositional logic follows from UIP of its modal companion.
Kurahashi \cite{Kur24} extended Visser's result and proved ULIP of $\mathbf{Int}$ and $\mathbf{KC}$, intermediate logics without local tabularity, by proving that of their respective modal companions $\mathbf{Grz}$ and $\mathbf{Grz.2}$.
Kogure and Kurahashi \cite{KK23} also proved UIP of a bimodal logic by constructing an embedding of it into a unimodal logic with UIP.
In this paper, we shall adopt a similar approach to these studies; we will establish an embedding of the logics $\NAmn$, $\NpAmn$, and $\NRAmn$ into $\CL$,
with which we can prove ULIP of them by reducing it to that of $\CL$.
In fact, we will introduce in Section \ref{sec:p18n} a general method, called \emph{propositionalization}, that enables one to transfer ULIP of a weaker logic to a stronger one, given an embedding with certain properties is established.
Our method is syntactic, as long as the embedding is obtained syntactically and ULIP of the weaker logic is proven syntactically.

We note that Gabbay and Maksimova \cite[Chapter 14]{GM09} also discuss a model-theoretic translation of a logic into a first-order theory based on its relational semantics.

\def\fCenter{\ \To\ }

\section{The sequent calculi for $\NAmn$, $\NpAmn$, and $\NRAmn$}\label{sec:sequent}

In this section, we introduce sequent calculi for logics $\NAmn$, $\NpAmn$, and $\NRAmn$, and prove that they enjoy cut elimination.

\begin{definition} \label{def:GNxAmn}
  Let $m, n \in \mathbb{N}$.
  The calculus $\GNAmn$ has the following initial sequents and rules.

  \subparagraph*{Initial sequents}
  Initial sequents have the form of $\varphi \To \varphi$ or $\bot \To$.

  \subparagraph*{Rules}
  Every rule in $\LK$ is a rule of $\GNAmn$.
  In addition to that, it has the following logical rules depending on $m$ and $n$:

  \begin{center}
    \Axiom$\fCenter \varphi$
    \RightLabel{($\mathrm{nec}$)}
    \UnaryInf$\fCenter \Box \varphi$
    \DisplayProof
  \end{center}
  \begin{center}
    \Axiom$\Box^m \varphi, \Box^n \varphi, \Gamma \fCenter \Delta$
    \RightLabel{(accL, when $n > m$)}
    \UnaryInf$\Box^n \varphi, \Gamma \fCenter \Delta$
    \DisplayProof
  \end{center}
  \begin{center}
    \Axiom$\Gamma \fCenter \Delta, \Box^m \varphi, \Box^n \varphi$
    \RightLabel{(accR, when $m > n$)}
    \UnaryInf$\Gamma \fCenter \Delta, \Box^m \varphi$
    \DisplayProof
  \end{center}

  $\GNpAmn$ is obtained from $\GNAmn$ by adding the following rule:
  \begin{center}
    \Axiom$\Box \varphi \fCenter$
    \RightLabel{(rosbox, when $m = 0$ and $n \ge 2$)}
    \UnaryInf$\Box \Box \varphi \fCenter$
    \DisplayProof
  \end{center}

  $\GNRAmn$ is obtained from $\GNAmn$ by adding the following rule:
  \begin{center}
    \Axiom$\varphi \fCenter$
    \RightLabel{(ros)}
    \UnaryInf$\Box \varphi \fCenter$
    \DisplayProof
  \end{center}
\end{definition}

For the sake of brevity, let $\NxAmn \in \set{ \NAmn,\, \NpAmn,\, \NRAmn }$,
and let $\GNxAmn$ designate the respective sequent calculus, $\GNAmn$, $\GNpAmn$, or $\GNRAmn$.

\begin{proposition} \label{prop:gnxamn-equiv}
  The following are equivalent:
  \begin{itemize}
    \item A sequent $\Gamma \To \Delta$ is provable in $\GNxAmn$.
    \item $\NxAmn \vdash \bigwedge \Gamma \to \bigvee \Delta$.
  \end{itemize}
  \begin{proof}
    Easy induction on the length of the proof.
  \end{proof}
\end{proposition}

\begin{corollary} \label{cor:gnxamn-consis}
  An empty sequent $\To$ is not provable in $\GNxAmn$.
  \begin{proof}
    Follows from the soundness of it (\cite[Theorem 4.4]{KS24}, Proposition \ref{prop:nramn-sound}).
  \end{proof}
\end{corollary}

\begin{theorem} \label{thm:gnxamn-cut-elim}
  $\GNxAmn$ enjoys cut elimination.
  \begin{proof}
    Take any proof of any sequent $\Gamma \To \Delta$ in $\GNxAmn$, and let
    \begin{center}
      \AxiomC{$\stackrel{\text{(cut-free)}}{\vdots}$}
      \RightLabel{\scriptsize{(LHS)}}
      \UnaryInfC{$\Gamma_1 \To \Delta_1, \varphi$}
      \AxiomC{$\stackrel{\text{(cut-free)}}{\vdots}$}
      \RightLabel{\scriptsize{(RHS)}}
      \UnaryInfC{$\varphi, \Gamma_2 \To \Delta_2$}
      \RightLabel{\scriptsize{(cut)}}
      \BinaryInfC{$\Gamma_1, \Gamma_2 \To \Delta_1, \Delta_2$}
      \DisplayProof
    \end{center}
    be one of the uppermost occurrences of the cut rule in the proof of $\Gamma \To \Delta$.
    Here, proofs marked `(cut-free)' refer to those that do not use the cut rule.
    We use an induction on the construction of the cut formula $\varphi$,
    and then an induction on the length of the proof.

    For the sake of brevity, we shall refer to the rule used to obtain $\Gamma_1 \To \Delta_1, \varphi$ as the \emph{LHS rule},
    and that of $\varphi, \Gamma_2 \To \Delta_2$ as the \emph{RHS rule}, given that they are not initial sequents.
    The following cases are proven similarly to that of $\LK$:
    \begin{itemize}
      \item At least one of $\Gamma_1 \To \Delta_1, \varphi$ or $\varphi, \Gamma_2 \To \Delta_2$ is an initial sequent.
      \item At least one of the LHS and the RHS is a structural rule.
      \item $\varphi$ is not the principal formula of at least one of the LHS and the RHS.
      \item Both the LHS and the RHS are the rules in $\LK$, and $\varphi$ is the principal formula of both.
    \end{itemize}
    Since every new rule in $\GNxAmn$ has a principal formula of the form $\Box \psi$,
    it suffices to consider the cases when both the LHS and the RHS are the new rules in $\GNxAmn$.
    Here, it is impossible to cut with the LHS being accR and the RHS being accL since they would not be the rules of $\GNxAmn$ at the same time.
    It is also impossible to cut with the LHS being accR and the RHS being rosbox by the same reason.
    We shall consider the following remaining cases.

    \paragraph*{Case 1:} Suppose that $\GNxAmn = \GNpAmn$, $m = 0$, $n \ge 2$, and the occurrence of the cut rule has the following form:
    \begin{center}
      \AxiomC{$\stackrel{\text{(cut-free)}}{\vdots}$}
      \noLine
      \UnaryInf$\fCenter \Box \psi$
      \RightLabel{\scriptsize{(nec)}}
      \UnaryInf$\fCenter \Box \Box \psi$
      \AxiomC{$\stackrel{\text{(cut-free)}}{\vdots}$}
      \noLine
      \UnaryInf$\Box \psi \fCenter$
      \RightLabel{\scriptsize{(rosbox)}}
      \UnaryInf$\Box \Box \psi \fCenter$
      \RightLabel{\scriptsize{(cut)}}
      \BinaryInfC{$\vphantom{\psi} \To$}
      \DisplayProof
    \end{center}
    This contradicts Corollary \ref{cor:gnxamn-consis}, so this case is impossible.

    \paragraph*{Case 2:} Suppose that $\GNxAmn = \GNRAmn$ and the occurrence of the cut rule has the following form:
    \begin{center}
      \AxiomC{$\stackrel{\text{(cut-free)}}{\vdots}$}
      \noLine
      \UnaryInf$\fCenter \psi$
      \RightLabel{\scriptsize{(nec)}}
      \UnaryInf$\fCenter \Box \psi$
      \AxiomC{$\stackrel{\text{(cut-free)}}{\vdots}$}
      \noLine
      \UnaryInf$\psi \fCenter$
      \RightLabel{\scriptsize{(ros)}}
      \UnaryInf$\Box \psi \fCenter$
      \RightLabel{\scriptsize{(cut)}}
      \BinaryInfC{$\vphantom{\psi} \To$}
      \DisplayProof
    \end{center}
    This also contradicts Corollary \ref{cor:gnxamn-consis}, so this case is also impossible.

    \paragraph*{Case 3:} Suppose that $n > m$ and the occurrence of the cut rule has the following form:
    \begin{center}
      \AxiomC{$\stackrel{\text{(cut-free)}}{\vdots}$}
      \noLine
      \UnaryInf$\fCenter \Box^{n-1} \psi$
      \RightLabel{\scriptsize{(nec)}}
      \UnaryInf$\fCenter \Box^n \psi$

      \AxiomC{$\stackrel{\text{(cut-free)}}{\vdots}$}
      \noLine
      \UnaryInf$\Box^m \psi, \Box^n \psi, \Gamma_2 \fCenter \Delta_2$
      \RightLabel{\scriptsize{(accL)}}
      \UnaryInf$\Box^n \psi, \Gamma_2 \fCenter \Delta_2$

      \RightLabel{\scriptsize{(cut)}}
      \BinaryInfC{$\Gamma_2 \To \Delta_2$}
      \DisplayProof
    \end{center}
    We shall refer to this proof as the \emph{original proof}.
    We first consider the following proof:
    \begin{center}
      \AxiomC{$\stackrel{\text{(cut-free)}}{\vdots}$}
      \noLine
      \UnaryInf$\fCenter \Box^{n-1} \psi$
      \RightLabel{\scriptsize{(nec)}}
      \UnaryInf$\fCenter \Box^n \psi$

      \AxiomC{$\stackrel{\text{(cut-free)}}{\vdots}$}
      \noLine
      \UnaryInf$\Box^m \psi, \Box^n \psi, \Gamma_2 \fCenter \Delta_2$

      \RightLabel{\scriptsize{(cut)}}
      \BinaryInfC{$\Box^m \psi, \Gamma_2 \To \Delta_2$}
      \DisplayProof
    \end{center}
    This proof is shorter than the original proof,
    so $\Box^m \psi, \Gamma_2 \To \Delta_2$ has a cut-free proof
    by the induction hypothesis on the length of the proof.

    Recall that $\To \Box^n \psi$ has a cut-free proof,
    then we can show that $\To \Box^m \psi$ also has a cut-free proof
    by Corollary \ref{cor:gnxamn-consis} and an easy induction on the length of the proof.
    Now we consider the following proof:
    \begin{center}
      \AxiomC{$\stackrel{\text{(cut-free)}}{\vdots}$}
      \noLine
      \UnaryInfC{$\fCenter \Box^m \psi$}

      \AxiomC{$\stackrel{\text{(cut-free)}}{\vdots}$}
      \noLine
      \UnaryInfC{$\Box^m \psi, \Gamma_2 \To \Delta_2$}

      \RightLabel{\scriptsize{(cut)}}
      \BinaryInfC{$\Gamma_2 \To \Delta_2$}
      \DisplayProof
    \end{center}
    The cut formula $\Box^m \psi$ of this proof is a proper subformula of the cut formula $\Box^n \psi$ of the original proof,
    so $\Gamma_2 \To \Delta_2$ has a cut-free proof
    by the induction hypothesis on the construction of the cut formula.

    \paragraph*{Case 4:} Suppose that $\GNxAmn = \GNRAmn$, $m > n$, and the occurrence of the cut rule has the following form:
    \begin{center}
      \AxiomC{$\stackrel{\text{(cut-free)}}{\vdots}$}
      \noLine
      \UnaryInf$\Gamma_2 \fCenter \Delta_2, \Box^m \varphi, \Box^n \varphi$
      \RightLabel{\scriptsize{(accR)}}
      \UnaryInf$\Gamma_2 \fCenter \Delta_2, \Box^m \varphi$

      \AxiomC{$\stackrel{\text{(cut-free)}}{\vdots}$}
      \noLine
      \UnaryInf$\Box^{m-1} \psi \fCenter$
      \RightLabel{\scriptsize{(ros)}}
      \UnaryInf$\Box^m \psi \fCenter$

      \RightLabel{\scriptsize{(cut)}}
      \BinaryInfC{$\Gamma_2 \To \Delta_2$}
      \DisplayProof
    \end{center}
    This case is proven similarly to the above.
  \end{proof}
\end{theorem}

\begin{remark}
  The reader may be wondering why we did not just use an initial sequent $\Box^n \varphi \To \Box^m \varphi$ to represent the $\Amn$ axiom, and why we introduced the accL and accR rules instead.
  Suppose $m > 0$ and $n = 0$, and consider the sequent calculus obtained from $\LK$ by adding the nec rule and the said initial sequent.
  This would permit the following cut, which cannot be eliminated:
  \begin{center}
    \def\fCenter{\ \To\ }

    \Axiom$\varphi_1 \fCenter \varphi_1$
    \RightLabel{\scriptsize{(\SequentRuleName{\lor}{R})}}
    \UnaryInf$\varphi_1 \fCenter \varphi_1 \lor \varphi_2$

    \Axiom$\varphi_1 \lor \varphi_2 \fCenter \Box^m (\varphi_1 \lor \varphi_2)$

    \RightLabel{\scriptsize{(cut)}}
    \BinaryInf$\varphi_1 \fCenter \Box^m (\varphi_1 \lor \varphi_2)$
    \DisplayProof
  \end{center}
  The same problem happens for the case when $m = 0$ and $n > 0$.
\end{remark}

Now we shall prove that $\NxAmn$ enjoys LIP by Maehara's method.

\begin{theorem}[Maehara's method] \label{thm:gnxamn-maehara}
  Take any sequent $\Gamma \To \Delta$ that is provable in $\GNxAmn$,
  then for any partition $(\Gamma_1 \To \Delta_1; \Gamma_2 \To \Delta_2)$ of $\Gamma \To \Delta$,
  there is $\chi \in \MF$ such that:
  \begin{enumerate}[label=(\alph*)]
    \item $\Gamma_1 \To \Delta_1, \chi$ is provable;
    \item $\chi, \Gamma_2 \To \Delta_2$ is provable;
    \item $\VarPos(\chi) \subseteq (\VarPos(\Gamma_1) \cup \VarNeg(\Delta_1)) \cap (\VarNeg(\Gamma_2) \cup \VarPos(\Delta_2))$;
    \item $\VarNeg(\chi) \subseteq (\VarNeg(\Gamma_1) \cup \VarPos(\Delta_1)) \cap (\VarPos(\Gamma_2) \cup \VarNeg(\Delta_2))$.
  \end{enumerate}

  \begin{proof}
    By Theorem \ref{thm:gnxamn-cut-elim}, there is a cut-free proof of $\Gamma \To \Delta$.
    We use an induction on the length of the cut-free proof.

    Take any partition $(\Gamma_1 \To \Delta_1; \Gamma_2 \To \Delta_2)$ of $\Gamma \To \Delta$.
    If $\Gamma \To \Delta$ is an initial sequent or is obtained by a rule of $\LK$, then
    it is proven by the same way one would apply Maehara's method to $\LK$.
    So we shall only consider the cases where $\Gamma \To \Delta$ is obtained by one of the new rules in $\GNxAmn$.

    If $(\Gamma \To \Delta) = (\To \Box \varphi)$ is obtained by $\frac{\To \varphi}{\To \Box \varphi}$ (nec), then
    either $\Box \varphi \in \Delta_1$ or $\Box \varphi \in \Delta_2$.
    We let $\chi = \bot$ or $\chi = \top$ respectively,
    and we can easily verify that $\chi$ satisfies (a) to (d) in both cases.

    If $\GNxAmn = \GNpAmn$, $m = 0$, $n \ge 2$, and $(\Gamma \To \Delta) = (\Box \Box \varphi \To)$ is obtained by $\frac{\Box \varphi \To}{\Box \Box \varphi \To}$ (rosbox),
    then either $\Box \Box \varphi \in \Gamma_1$ or $\Box \Box \varphi \in \Gamma_2$.
    We let $\chi = \bot$ or $\chi = \top$ respectively,
    and we can easily verify that $\chi$ satisfies (a) to (d) in both cases.

    If $\GNxAmn = \GNRAmn$ and $(\Gamma \To \Delta) = (\Box \varphi \To)$ is obtained by $\frac{\varphi \To}{\Box \varphi \To}$ (ros),
    then it is proven similarly to the above case.

    If $n > m$ and $(\Gamma \To \Delta) = (\Box^n \varphi, \Gamma' \To \Delta)$ is obtained by $\frac{\Box^m \varphi, \Box^n \varphi, \Gamma' \To \Delta}{\Box^n \varphi, \Gamma' \To \Delta}$ (accL),
    then either $\Box^n \varphi \in \Gamma_1$ or $\Box^n \varphi \in \Gamma_2$.

    \begin{itemize}
      \item If $\Box^n \varphi \in \Gamma_1$, then let $\Gamma_1' = \Gamma_1 \setminus \set{\Box^n \varphi}$.
            Here, $(\Box^m \varphi, \Box^n \varphi, \Gamma_1' \To \Delta_1; \Gamma_2 \To \Delta_2)$ is a partition of $\Box^m \varphi, \Box^n \varphi, \Gamma' \To \Delta$.
            Then by the induction hypothesis, there is $\chi \in \MF$ such that:
            \begin{itemize}
              \item $\Box^m \varphi, \Box^n \varphi, \Gamma_1' \To \Delta_1, \chi$ is provable, which implies (a) by applying the accL rule;
              \item $\chi, \Gamma_2 \To \Delta_2$ is provable, which is (b);
              \item $\VarPos(\chi) \subseteq (\VarPos(\set{\Box^m \varphi, \Box^n \varphi} \cup \Gamma_1') \cup \VarNeg(\Delta_1)) \cap (\VarNeg(\Gamma_2) \cup \VarPos(\Delta_2))$,
                    which implies (c) along with $\VarPos(\set{\Box^m \varphi, \Box^n \varphi}) = \VarPos(\Box^n \varphi)$;
              \item $\VarNeg(\chi) \subseteq (\VarNeg(\set{\Box^m \varphi, \Box^n \varphi} \cup \Gamma_1') \cup \VarPos(\Delta_1)) \cap (\VarPos(\Gamma_2) \cup \VarNeg(\Delta_2))$,
                    which implies (d) along with $\VarNeg(\set{\Box^m \varphi, \Box^n \varphi}) = \VarNeg(\Box^n \varphi)$.
            \end{itemize}
      \item If $\Box^n \varphi \in \Gamma_2$, then let $\Gamma_2' = \Gamma_2 \setminus \set{\Box^n \varphi}$
            and consider a partition $(\Gamma_1 \To \Delta_1; \Box^m \varphi, \Box^n \varphi, \Gamma_2' \To \Delta_2)$ of $\Box^m \varphi, \Box^n \varphi, \Gamma' \To \Delta$,
            then the rest of the proof is similar to the previous case.
    \end{itemize}

    If $m > n$ and $(\Gamma \To \Delta) = (\Gamma \To \Delta', \Box^m \varphi)$ is obtained by $\frac{\Gamma \To \Delta', \Box^m \varphi, \Box^n \varphi}{\Gamma \To \Delta', \Box^m \varphi}$ (accR),
    then either $\Box^m \varphi \in \Delta_1$ or $\Box^m \varphi \in \Delta_2$.
    \begin{itemize}
      \item If $\Box^m \varphi \in \Delta_1$, then let $\Delta_1' = \Delta_1 \setminus \set{\Box^m \varphi}$
            and consider a partition $(\Gamma_1 \To \Delta_1', \Box^m \varphi, \Box^n \varphi; \Gamma_2 \To \Delta_2)$ of $\Gamma \To \Delta', \Box^m \varphi, \Box^n \varphi$.
      \item If $\Box^m \varphi \in \Delta_2$, then let $\Delta_2' = \Delta_2 \setminus \set{\Box^m \varphi}$
            and consider a partition $(\Gamma_1 \To \Delta_1; \Gamma_2 \To \Delta_2', \Box^m \varphi, \Box^n \varphi)$ of $\Gamma \To \Delta', \Box^m \varphi, \Box^n \varphi$.
    \end{itemize}
    In any case, the rest of the proof is similar.
  \end{proof}
\end{theorem}

\begin{corollary} \label{cor:nramn-lip}
  $\NxAmn$ enjoys LIP and, consequently, enjoys CIP.
\end{corollary}

\section{Propositionalization} \label{sec:p18n}

In this section, we introduce a general method, called \emph{propositionalization},
that, with certain conditions satisfied, enables us to prove ULIP of a logic by reducing it to that of some weaker logic.

Let $L, M$ be logics with languages $\Lang_L$ and $\Lang_M$, respectively, where $\PF \subseteq \Lang_L \subseteq \Lang_M \subseteq \MF$,
such that for every $\varphi \in \Lang_L$ and every substitution $\sigma: \PropVar \to \Lang_M$,
if $L \vdash \varphi$, then $M \vdash \sigma(\varphi)$.
We moreover assume that $M$ is closed under syllogism:
if $M \vdash \varphi \to \psi$ and $M \vdash \psi \to \chi$, then $M \vdash \varphi \to \chi$.

\begin{definition} \label{def:extend-lang}
  Let $\PropVar' = \PropVar \cup \Set{ p_\varphi; p_\varphi \notin \PropVar, \varphi \in \Lang_M }$.
  We extend $\Lang_L$ to $\Lang_L'$ by extending its propositional variables to $\PropVar'$.
  We also define $\VarAny: \Lang_L' \to \PropVar'$ by extending $\VarAny: \Lang_L \to \PropVar$ appropriately.
  We denote by $L'$ the logic obtained from $L$ by replacing its language $\Lang_L$ with $\Lang_L'$.
\end{definition}

\begin{definition} \label{def:stdsub}
  We define a \emph{standard substitution} $\sigma^{\text{S}}: \PropVar' \to \Lang_M$ by:
  \begin{equation*}
    \sigma^{\text{S}}(p) =
    \left\{\begin{array}{ll}
              p       & \mathrel{\text{if}} p \in \PropVar, \\
              \varphi & \mathrel{\text{if}} p = p_\varphi \mathrel{\text{for some}} \varphi \in \Lang_M.
            \end{array}
    \right.
  \end{equation*}
  Here, it is clear that $L' \vdash \rho$ implies $M \vdash \sigma^{\text{S}}(\rho)$ for any $\rho \in \Lang_L'$.
\end{definition}

\begin{definition} \label{def:p18n}
  Let $\sharp, \flat: \Lang_M \to \Lang_L'$. We say the pair $(\sharp, \flat)$ is a \emph{propositionalization} of $M$ into $L$
  if for any $\varphi, \psi \in \Lang_M$, the following hold:
  \begin{enumerate}[label=(\arabic*)]
    \item $M \vdash \sigma^{\text{S}}(\varphi^\sharp) \to \varphi$ and $M \vdash \varphi \to \sigma^{\text{S}}(\varphi^\flat)$;
    \item $M \vdash \varphi \to \psi$ implies $L' \vdash \varphi^\flat \to \psi^\sharp$;
    \item For $(\bullet, \circ) \in \set{ (+, -), (-, +) }$ and $\natural \in \set{\sharp, \flat}$,
          \begin{enumerate}
            \item $p \in \VarAny(\varphi^\natural) \cap \PropVar$ implies $p \in \VarAny(\varphi)$;
            \item $p_\psi \in \VarAny(\varphi^\natural) \setminus \PropVar$ implies $\VarAny(\psi) \subseteq \VarPos(\varphi)$ and $\VarAnother(\psi) \subseteq \VarNeg(\varphi)$.
          \end{enumerate}
  \end{enumerate}
\end{definition}

\begin{theorem} \label{thm:p18n-ulip}
  If there is a propositionalization $(\sharp, \flat)$ of $M$ into $L$ and $L$ enjoys ULIP, then so does $M$.
  \begin{proof}
    Suppose that $L$ enjoys ULIP, then $L'$ does also.
    To prove ULIP of $M$, take any $\varphi \in \Lang_M$ and any finite $P^+, P^- \subseteq \PropVar$.

    Let us say $\psi \in \Lang_M$ is \emph{$+$-safe} if $P^+ \cap \VarPos(\psi) = P^- \cap \VarNeg(\psi) = \emptyset$,
    and is \emph{$-$-safe} if $P^+ \cap \VarNeg(\psi) = P^- \cap \VarPos(\psi) = \emptyset$.
    For $\bullet \in \set{+, -}$, we let:
    \begin{equation*}
      Q^\bullet = P^\bullet \cup \Set{ p_\psi \in \Var(\varphi^\flat); \psi \mathrel{\text{is not}} \bullet\text{-safe} }.
    \end{equation*}
    Then by ULIP of $L'$, there is $\chi' \in \Lang_L'$ such that
    $\VarAny(\chi') \subseteq \VarAny(\varphi^\flat) \setminus Q^\bullet$,
    $L' \vdash \varphi^\flat \to \chi'$,
    and $L' \vdash \chi' \to \psi'$ for any $\psi' \in \Lang_L'$ such that $L' \vdash \varphi^\flat \to \psi'$ and $\VarAny(\psi') \cap Q^\bullet = \emptyset$, for $\bullet \in \set{+, -}$.
    We let $\chi = \sigma^{\text{S}}(\chi')$ and we shall prove that $\chi$ satisfies the following conditions:
    (a) $\VarAny(\chi) \subseteq \VarAny(\varphi) \setminus P^\bullet$, for $\bullet \in \set{+, -}$,
    (b) $M \vdash \varphi \to \chi$, and
    (c) $M \vdash \chi \to \psi$ for any $\psi \in \Lang_M$ such that $M \vdash \varphi \to \psi$ and $\VarAny(\psi) \cap P^\bullet = \emptyset$ for $\bullet \in \set{+, -}$.

    \subparagraph*{(a)}
    We first evaluate $\VarAny(\chi')$ for $\bullet \in \set{+, -}$:
    \begin{align*}
      \VarAny(\chi')
       & \subseteq \VarAny(\varphi^\flat) \setminus Q^\bullet                                                                              \\
       & \subseteq \left( \VarAny(\varphi) \cup \set*{ p_\rho; p_\rho \in \VarAny(\varphi^\flat) } \right) \setminus Q^\bullet             \\
       & = \left(\VarAny(\varphi) \setminus P^\bullet\right) \cup \set*{ p_\rho \in \VarAny(\varphi^\flat); \rho \mathrel{\text{is}} \bullet\text{-safe} }.
    \end{align*}

    We note that the second inclusion is by Definition \ref{def:p18n} (3a).
    Now we shall prove that $\VarAny(\chi) \subseteq \VarAny(\varphi) \setminus P^\bullet$
    for $(\bullet, \circ) \in \set{ (+, -), (-, +) }$.
    Take any $p \in \VarAny(\chi)$.
    If $p \in \VarAny(\chi')$, then $p \in \VarAny(\varphi) \setminus P^\bullet$ by the above.
    Otherwise, by $\chi = \sigma^{\text{S}}(\chi')$, there must be $p_\rho \in \Var(\chi')$ such that $p \in \Var(\rho)$.
    Since $p \in \VarAny(\chi)$, at least one of the following holds:
    (i)  $p \in \VarAny(\rho)$ and $p_\rho \in \VarPos(\chi')$, or
    (ii) $p \in \VarAnother(\rho)$ and $p_\rho \in \VarNeg(\chi')$.
    \begin{enumerate}[label=(\roman*)]
      \item If $p \in \VarAny(\rho)$ and $p_\rho \in \VarPos(\chi')$,
            then $p_\rho \in \VarPos(\varphi^\flat)$ and $\rho$ is $+$-safe by the above.
            By Definition \ref{def:p18n} (3b), $p_\rho \in \VarPos(\varphi^\flat)$ implies $\VarAny(\rho) \subseteq \VarAny(\varphi)$, so $p \in \VarAny(\varphi)$.
            Also, $\rho$ being $+$-safe implies $\VarAny(\rho) \cap P^\bullet = \emptyset$, so $p \notin P^\bullet$.
      \item If $p \in \VarAnother(\rho)$ and $p_\rho \in \VarNeg(\chi')$,
            then $p_\rho \in \VarNeg(\varphi^\flat)$ and $\rho$ is $-$-safe by the above.
            By Definition \ref{def:p18n} (3b), $p_\rho \in \VarNeg(\varphi^\flat)$ implies $\VarAnother(\rho) \subseteq \VarAny(\varphi)$, so $p \in \VarAny(\varphi)$.
            Also, $\rho$ being $-$-safe implies $\VarAnother(\rho) \cap P^\bullet = \emptyset$, so $p \notin P^\bullet$.
    \end{enumerate}

    \subparagraph*{(b)}
    $L' \vdash \varphi^\flat \to \chi'$ implies $M \vdash \sigma^{\text{S}}(\varphi^\flat) \to \chi$ by Definition \ref{def:stdsub}.
    Also, $M \vdash \varphi \to \sigma^{\text{S}}(\varphi^\flat)$ holds by Definition \ref{def:p18n} (1).
    Therefore, $M \vdash \varphi \to \chi$ holds by syllogism.

    \subparagraph*{(c)}
    Take any such $\psi \in \Lang_M$, then $M \vdash \varphi \to \psi$ and $\VarAny(\psi) \cap P^\bullet = \emptyset$ hold for $\bullet \in \set{+, -}$.
    The former implies $L' \vdash \varphi^\flat \to \psi^\sharp$ by Definition \ref{def:p18n} (2).
    With the latter in mind, we first evaluate $\VarAny(\psi^\sharp) \cap Q^\bullet$ for $\bullet \in \set{+, -}$:
    \begin{align*}
      \VarAny(\psi^\sharp) \cap Q^\bullet
       & \subseteq \left(\VarAny(\psi) \cup \set*{ p_\rho; p_\rho \in \VarAny(\psi^\sharp) }\right) \cap Q^\bullet                                     \\
       & = \left(\VarAny(\psi) \cap P^\bullet\right) \cup \set*{ p_\rho; p_\rho \in \VarAny(\psi^\sharp) \cap \VarAny(\varphi^\flat), \rho \mathrel{\text{is not}} \bullet\text{-safe} } \\
       & \subseteq \set*{ p_\rho \in \VarAny(\psi^\sharp); \rho \mathrel{\text{is not}} \bullet\text{-safe} }.
    \end{align*}

    We note that the first inclusion is by Definition \ref{def:p18n} (3a).
    We let $R^\bullet = \set*{ p_\rho \in \VarAny(\psi^\sharp); \rho \mathrel{\text{is not}} \bullet\text{-safe} }$, and we shall show that $R^\bullet = \emptyset$ for $(\bullet, \circ) \in \set{ (+, -), (-, +) }$.
    Take any $p_\rho \in R^\bullet$, then $p_\rho \in \VarAny(\psi^\sharp)$ implies that $\VarAny(\rho) \subseteq \VarPos(\psi)$ and $\VarAnother(\rho) \subseteq \VarNeg(\psi)$.
    Here, $P^+ \cap \VarAny(\rho) \subseteq P^+ \cap \VarPos(\psi) = \emptyset$ and $P^- \cap \VarAnother(\rho) \subseteq P^- \cap \VarNeg(\psi) = \emptyset$,
    so $\rho$ is $\bullet$-safe, which is a contradiction. Therefore, $R^\bullet = \emptyset$.
    Now that $\VarAny(\psi^\sharp) \cap Q^\bullet \subseteq R^\bullet = \emptyset$ for $\bullet \in \set{+, -}$,
    then $L' \vdash \chi' \to \psi^\sharp$ by $\chi'$ being the post-interpolant of $(\varphi^\flat, Q^+, Q^-)$,
    so $M \vdash \chi \to \sigma^{\text{S}}(\psi^\sharp)$.
    Also, $M \vdash \sigma^{\text{S}}(\psi^\sharp) \to \psi$ holds by Definition \ref{def:p18n} (1).
    Therefore, $M \vdash \chi \to \psi$ holds by syllogism.
  \end{proof}
\end{theorem}

\begin{example}
  We shall present a trivial example of propositionalization.
  Let $N: \PF \to \PF'$ be defined by $\varphi^N = \neg \neg \varphi$ for $\varphi \in \PF$, then 
  we can verify that $(N, N)$ is a propositionalization of $\CL$ into intuitionistic propositional logic $\mathbf{Int}$:
  \begin{enumerate}
    \item It is easy to see that $\sigma^{\text{S}}(\varphi^N) = \neg \neg \varphi$ since $p_\psi \notin \VarAny(\varphi^N)$ for any $\psi \in \PF$. Here, both $\CL \vdash \neg \neg \varphi \to \varphi$ and $\CL \vdash \varphi \to \neg \neg \varphi$ hold.
    \item Suppose that $\CL \vdash \varphi \to \psi$, then $\mathbf{Int}' \vdash \neg \neg (\varphi \to \psi)$ by Glivenko's theorem,
          which implies $\mathbf{Int}' \vdash (\neg \neg \varphi) \to (\neg \neg \psi)$,
          so $\mathbf{Int}' \vdash \varphi^N \to \psi^N$.
    \item $\VarAny(\varphi^N) = \VarAny(\varphi) \subseteq \PropVar$ trivially implies both conditions.
  \end{enumerate}
  So by Theorem \ref{thm:p18n-ulip}, we can see that ULIP of $\mathbf{Int}$ implies that of $\CL$.
\end{example}

\begin{example}
  The boxdot translation $(\Box \varphi)^\boxdot = \varphi^\boxdot \land \Box \varphi^\boxdot$ (see, e.g., \cite{Jer16}) is another trivial example of propositionalization.
  Let $L, M$ be modal logics where $L \subseteq M$, $\mathbf{KT} \subseteq M$, and $L$ is a boxdot companion of $M$ (i.e.\ $M \vdash \varphi$ iff $L \vdash \varphi^\boxdot$),
  then we can easily check that $(\boxdot, \boxdot)$ is a propositionalization of $M$ into $L$.
  So by Theorem \ref{thm:p18n-ulip}, ULIP of $L$ implies ULIP of $M$.
  For example, ULIP of $\K$ implies ULIP of $\mathbf{KT}$, and the failure of ULIP in $\mathbf{S4}$ \cite{GZ95} implies that for any boxdot companion $L$ of $\mathbf{S4}$, if $L \subseteq \mathbf{S4}$, then $L$ also lacks ULIP.
\end{example}

In particular, we obtain the following corollary of Theorem \ref{thm:p18n-ulip}.

\begin{corollary}
  Let $L, M$ be logics with the same language where $L \subseteq M$ and $M$ is closed under syllogism.
  If there is a pair of translations $(\sharp, \flat)$ such that
  1) $M \vdash \varphi^\sharp \to \varphi$ and $M \vdash \varphi \to \varphi^\flat$,
  2) $M \vdash \varphi \to \psi$ implies $L \vdash \varphi^\flat \to \psi^\sharp$,
  and 3) $\VarAny(\varphi^\natural) \subseteq \VarAny(\varphi)$ for $\bullet \in \set{+, -}, \natural \in \set{\sharp, \flat}$,
  and $L$ enjoys ULIP, then so does $M$.
\end{corollary}

In this paper, we are particularly interested in the case when $L = \CL$ and $M \in \Set{ \NAmn,\, \NpAmn,\, \NRAmn }$.
However, we believe that this method can also be used to analyze ULIP of intuitionistic variants of $\N$ and its extensions
by reducing it to that of $\mathbf{Int}$. Moreover, given that the definition of $\VarAny$ is extended properly for bimodal language,
it would enable one to analyze ULIP of a bimodal logic by propositionalizing it to a unimodal logic, generalizing the work of Kogure and Kurahashi \cite{KK23}.

\def\fCenter{\ \To\ }

\section{ULIP of $\NAmn$, $\NpAmn$, and $\NRAmn$} \label{sec:ulip}

In this section, we prove ULIP of $\NAmn$, $\NpAmn$, and $\NRAmn$ by constructing propositionalizations into $\CL$.

As in Section \ref{sec:sequent}, let $\NxAmn \in \set{ \NAmn,\, \NpAmn,\, \NRAmn }$,
and let $\GNxAmn$ designate the respective sequent calculus, $\GNAmn$, $\GNpAmn$, or $\GNRAmn$.
We shall briefly check that $\NxAmn$ and $\CL$ satisfy the conditions we mentioned at the beginning of Section \ref{sec:p18n}.
\begin{itemize}
  \item $\PF \subseteq \MF$.
  \item For every $\varphi \in \PF$ and every substitution $\sigma: \PropVar \to \MF$, if $\CL \vdash \varphi$, then $\NxAmn \vdash \sigma(\varphi)$.
  \item Syllogism holds in $\NxAmn$.
\end{itemize}
Then, $\PF'$, $\CL'$, and $\VarAny: \PF' \to \PropVar'$ are obtained by Definition \ref{def:extend-lang}.

Now, we shall first create a pair of translations $(\sharp, \flat)$ in such a way that we can \emph{emulate} every rule of $\GNxAmn$ in $\LK$.
Given a modal sequent $\Gamma \To \Delta$, we translate it to a propositional one $\Gamma^\flat \To \Delta^\sharp$ by
replacing every subformula $\Box \psi$ by $p_{\Box \psi} \in \PropVar'$ with extra information to emulate the accL and accR rules.
The nec and ros(box) rules would not be able to emulate in $\LK$
since their principal formulae would be the fresh variable that is not present in the upper sequent, so we replace it with $\To \top$ and $\bot \To$ respectively.
The following translations do exactly the above.

\begin{definition} \label{def:sharpflat}
  For $\varphi \in \MF$, we define $\varphi^\sharp, \varphi^\flat \in \PF'$ inductively:
  \begin{itemize}
    \item $\bot^\natural = \bot$, $p^\natural = p$, and $(\psi_1 \circledcirc \psi_2)^\natural = \psi_1^\natural \circledcirc \psi_2^\natural$, for $\natural \in \set{\sharp, \flat}$ and $\circledcirc \in \set{\land, \lor}$.
    \item $(\psi_1 \to \psi_2)^\sharp = \psi_1^\flat \to \psi_2^\sharp$ and $(\psi_1 \to \psi_2)^\flat = \psi_1^\sharp \to \psi_2^\flat$.
    \item For any $\Box^k \psi$ such that $k \ge 1$ and $\psi \ne \Box \chi$ for any $\chi \in \MF$:
          \begin{align*}
            (\Box^k \psi)^\sharp & =
            \left\{\begin{array}{ll}
                     \top                                            & \tif \NxAmn \vdash \Box^{k-1} \psi, \\
                     p_{\Box^k \psi} \lor (\Box^{k-m+n} \psi)^\sharp & \tif \NxAmn \nvdash \Box^{k-1} \psi \tand k \ge m > n,                   \\
                     p_{\Box^k \psi}                                 & \mathrel{\text{otherwise.}}
                   \end{array}
            \right.                  \\
            (\Box^k \psi)^\flat  & =
            \left\{\begin{array}{ll}
                     \bot                                            & \tif C(k, \psi),                                                    \\
                     p_{\Box^k \psi} \land (\Box^{k-n+m} \psi)^\flat & \tif\ \tnot C(k, \psi) \tand k \ge n > m,                                                    \\
                     p_{\Box^k \psi}                                 & \mathrel{\text{otherwise.}}
                   \end{array}
            \right.
          \end{align*}
  \end{itemize}
  Here, the condition $C(k, \psi)$ is defined by:
  \begin{itemize}
    \item If $\NxAmn = \NpAmn$, $m = 0$, $n \ge 2$, $k \ge 2$, and $\NpAmn \vdash \neg \Box^{k-1} \psi$, then $C(k, \psi)$ holds;
    \item If $\NxAmn = \NRAmn$ and $\NRAmn \vdash \neg \Box^{k-1} \psi$, then $C(k, \psi)$ holds;
    \item Otherwise, $C(k, \psi)$ does not hold.
  \end{itemize}
  For $\Gamma \subseteq \MF$ and $\natural \in \set{\sharp, \flat}$, we also let $\Gamma^\natural = \set{ \varphi^\natural; \varphi \in \Gamma }$.
\end{definition}

We shall then show that $(\sharp, \flat)$ is indeed a propositionalization of $\NxAmn$ into $\CL$,
proving ULIP of $\NxAmn$ as a consequence.

\begin{proposition}[3.\ of Def.\ \ref{def:p18n}] \label{prop:sharpflat-variables}
  For $(\bullet, \circ) \in \set{ (+, -), (-, +) }$ and $\natural \in \set{\sharp, \flat}$, the following hold:
  \begin{itemize}
    \item $p \in \VarAny(\varphi^\natural) \cap \PropVar$ implies $p \in \VarAny(\varphi)$;
    \item $p_\psi \in \VarAny(\varphi^\natural) \setminus \PropVar$ implies $\VarAny(\psi) \subseteq \VarPos(\varphi)$ and $\VarAnother(\psi) \subseteq \VarNeg(\varphi)$.
  \end{itemize}
  \begin{proof}
    Easy induction on the construction of $\varphi$.
  \end{proof}
\end{proposition}

\begin{proposition} \label{prop:flat-dominates-sharp}
  $\CL' \vdash \varphi^\flat \to \varphi^\sharp$ for any $\varphi \in \MF$.
  \begin{proof}
    We use an induction on the construction of $\varphi$.
    It holds trivially when $\varphi = \bot$ or $\varphi = p$ for $p \in \PropVar$.
    If $\varphi = \psi_1 \circledcirc \psi_2$, for $\circledcirc \in \set{ \land, \lor, \to }$, then it is easily proven by the induction hypothesis.
    Now suppose that $\varphi = \Box^k \psi$ ($k \ge 1$, $\psi \ne \Box \chi$).
    If $(\Box^k \psi)^\sharp = \top$ or $(\Box^k \psi)^\flat = \bot$,
    then $\CL' \vdash (\Box^k \psi)^\flat \to \top$ or $\CL' \vdash \bot \to (\Box^k \psi)^\sharp$ trivially hold, respectively.
    Suppose not, then we distinguish the following cases:
    \begin{itemize}
      \item If $(\Box^k \psi)^\sharp = p_{\Box^k \psi} \lor (\Box^{k-m+n} \psi)^\sharp$ and $(\Box^k \psi)^\flat = p_{\Box^k \psi}$ ($k \ge m > n$),
            then $\CL' \vdash p_{\Box^k \psi} \to (p_{\Box^k \psi} \lor (\Box^{k-m+n} \psi)^\sharp)$ holds.
      \item If $(\Box^k \psi)^\flat = p_{\Box^k \psi} \land (\Box^{k-n+m} \psi)^\flat$ and $(\Box^k \psi)^\sharp = p_{\Box^k \psi}$ ($k \ge n > m$),
            then $\CL' \vdash (p_{\Box^k \psi} \land (\Box^{k-n+m} \psi)^\flat) \to p_{\Box^k \psi}$ holds.
      \item If $(\Box^k \psi)^\sharp = (\Box^k \psi)^\flat = p_{\Box^k \psi}$,
            then $\CL' \vdash p_{\Box^k \psi} \to p_{\Box^k \psi}$ holds. \qedhere
    \end{itemize}
  \end{proof}
\end{proposition}

\begin{proposition}[1.\ of Def.\ \ref{def:p18n}] \label{prop:equality}
  $\NxAmn \vdash \sigma^{\text{S}}(\varphi^\sharp) \to \varphi$ and $\NxAmn \vdash \varphi \to \sigma^{\text{S}}(\varphi^\flat)$ for any $\varphi \in \MF$.
  \begin{proof}
    We use an induction on the construction of $\varphi$ to prove both of the above simultaneously.
    
    Both hold trivially when $\varphi = \bot$ or $\varphi = p$ for $p \in \PropVar$.
    
    If $\varphi = \psi_1 \circledcirc \psi_2$, for $\circledcirc \in \set{ \land, \lor }$, then both are easily proven by the induction hypothesis.

    Suppose that $\varphi = \psi_1 \to \psi_2$. We shall only prove $\NxAmn \vdash \sigma^S(\varphi^\sharp) \to \varphi$,
    as $\NxAmn \vdash \varphi \to \sigma^S(\varphi^\flat)$ is proven dually.
    Here, $\sigma^S(\varphi^\sharp) = \sigma^S(\psi_1^\flat \to \psi_2^\sharp) = \sigma^S(\psi_1^\flat) \to \sigma^S(\psi_2^\sharp)$,
    then by the induction hypothesis, we have $\NxAmn \vdash \psi_1 \to \sigma^S(\psi_1^\flat)$ and $\NxAmn \vdash \sigma^S(\psi_2^\sharp) \to \psi_2$.
    Now we have $\NxAmn \vdash \sigma^S(\varphi^\sharp) \to \varphi$ by the following derivation in $\NxAmn$:
    \begin{center}
    \scalebox{0.8}{
      \AxiomC{$ $}
      \RightLabel{\scriptsize{I.H.}}
      \UnaryInfC{$\psi_1 \to \sigma^S(\psi_1^\flat)$}

      \AxiomC{$ $}
      \UnaryInfC{$p \to ((p \to q) \to q)$}
      \RightLabel{\scriptsize{Subst.}}
      \UnaryInfC{$\sigma^S(\psi_1^\flat) \to ((\sigma^S(\psi_1^\flat) \to \sigma^S(\psi_2^\sharp)) \to \sigma^S(\psi_2^\sharp))$}
      \UnaryInfC{$\sigma^S(\psi_1^\flat) \to (\sigma^S(\psi_1^\flat \to \psi_2^\sharp) \to \sigma^S(\psi_2^\sharp))$}

      \AxiomC{$ $}
      \RightLabel{\scriptsize{I.H.}}
      \UnaryInfC{$\sigma^S(\psi_2^\sharp) \to \psi_2$}

      \RightLabel{\scriptsize{Syll.}}
      \TrinaryInfC{$\psi_1 \to (\sigma^S(\psi_1^\flat \to \psi_2^\sharp) \to \psi_2)$}
      \UnaryInfC{$\sigma^S(\psi_1^\flat \to \psi_2^\sharp) \to (\psi_1 \to \psi_2)$}
      \UnaryInfC{$\sigma^S(\varphi^\sharp) \to \varphi$}
      \DisplayProof
    }
    \end{center}

    Suppose that $\varphi = \Box^k \psi$ ($k \ge 1$, $\psi \ne \Box \chi$).
    We first prove $\NxAmn \vdash \sigma^{\text{S}}((\Box^k \psi)^\sharp) \to \Box^k \psi$.
    We distinguish the following cases.
    \begin{itemize}
      \item If $(\Box^k \psi)^\sharp = \top$,
            then $\sigma^{\text{S}}((\Box^k \psi)^{\sharp}) = \top$.
            Here, $\NxAmn \vdash \Box^k \psi$ by the definition of $\sharp$ and the $\textsc{Nec}$ rule,
            so $\NxAmn \vdash \top \to \Box^k \psi$ holds.
      \item If $k \ge m > n$ and $(\Box^k \psi)^\sharp = p_{\Box^k \psi} \lor (\Box^{k-m+n} \psi)^\sharp$,
            then $\sigma^{\text{S}}((\Box^k \psi)^\sharp) = \Box^k \psi \lor \sigma^{\text{S}}((\Box^{k-m+n} \psi)^\sharp)$.
            Then by the induction hypothesis, $\NxAmn \vdash \sigma^{\text{S}}((\Box^{k-m+n} \psi)^\sharp) \to \Box^{k-m+n} \psi$.
            Also, $\NxAmn \vdash \Box^{k-m+n} \psi \to \Box^k \psi$ by the $\Amn$ axiom.
            So $\NxAmn \vdash (\Box^k \psi \lor \sigma^{\text{S}}((\Box^{k-m+n} \psi)^\sharp)) \to \Box^k \psi$ holds.
      \item Otherwise, $(\Box^k \psi)^\sharp = p_{\Box^k \psi}$ and then $\sigma^{\text{S}}((\Box^k \psi)^\sharp) = \Box^k \psi$.
            Here, $\NxAmn \vdash \Box^k \psi \to \Box^k \psi$ trivially holds.
    \end{itemize}
    Now we prove $\NxAmn \vdash \Box^k \psi \to \sigma^{\text{S}}((\Box^k \psi)^\flat)$.
    We distinguish the following cases.
    \begin{itemize}
      \item If $\NxAmn = \NpAmn$, $m = 0$, $n \ge 2$, $k \ge 2$, and $(\Box^k \psi)^\flat = \bot$,
            then $\sigma^{\text{S}}((\Box^k \psi)^\flat) = \bot$.
            Here, $\NpAmn \vdash \neg \Box^k \psi$ by the definition of $\flat$ and the $\RosBox$ rule,
            so $\NpAmn \vdash \Box^k \psi \to \bot$ holds.
      \item If $\NxAmn = \NRAmn$ and $(\Box^k \psi)^\flat = \bot$,
            then $\sigma^{\text{S}}((\Box^k \psi)^\flat) = \bot$.
            Here, $\NRAmn \vdash \neg \Box^k \psi$ by the definition of $\flat$ and the $\Ros$ rule,
            so $\NRAmn \vdash \Box^k \psi \to \bot$ holds.
      \item If $k \ge n > m$ and $(\Box^k \psi)^\flat = p_{\Box^k \psi} \land (\Box^{k-n+m} \psi)^\flat$,
            then $\sigma^{\text{S}}((\Box^k \psi)^\flat) = \Box^k \psi \land \sigma^{\text{S}}((\Box^{k-n+m} \psi)^\flat)$.
            Then by the induction hypothesis, $\NxAmn \vdash \Box^{k-n+m} \psi \to \sigma^{\text{S}}((\Box^{k-n+m} \psi)^\flat)$.
            Also, $\NxAmn \vdash \Box^k \psi \to \Box^{k-n+m} \psi$ by the $\Amn$ axiom.
            So $\NxAmn \vdash  \Box^k \psi \to (\Box^k \psi \land \sigma^{\text{S}}((\Box^{k-n+m} \psi)^\flat))$ holds.
      \item Otherwise, $(\Box^k \psi)^\flat = p_{\Box^k \psi}$ and then $\sigma^{\text{S}}((\Box^k \psi)^\flat) = \Box^k \psi$.
            Here, $\NxAmn \vdash \Box^k \psi \to \Box^k \psi$ trivially holds. \qedhere
    \end{itemize}
  \end{proof}
\end{proposition}

Now we shall show that every rule in $\GNxAmn$ can indeed be emulated in $\LK$ using the $\sharp/\flat$-translations.

\begin{theorem}[Emulation] \label{thm:emulation}
  If a sequent $\Gamma \To \Delta$ is provable in $\GNxAmn$,
  then $\Gamma^\flat \To \Delta^\sharp$ is provable in $\LK$.
  \begin{proof}
    We use an induction on the length of the cut-free proof of $\Gamma \To \Delta$.

    If $(\Gamma \To \Delta) = (\bot \To)$ is an initial sequent,
    then $(\Gamma^\flat \To \Delta^\sharp) = (\bot \To)$ is also an initial sequent in $\LK$.

    If $(\Gamma \To \Delta) = (\varphi \To \varphi)$ is an initial sequent,
    then $(\Gamma^\flat \To \Delta^\sharp) = (\varphi^\flat \To \varphi^\sharp)$ is provable in $\LK$ by Proposition \ref{prop:flat-dominates-sharp}.

    If $\Gamma \To \Delta$ is obtained by a rule in $\LK$, then it is easily proven with the induction hypothesis.
    Here we demonstrate some of such cases.
    \begin{itemize}
      \item If $(\Gamma \To \Delta) = (\Gamma_1, \Gamma_2, \varphi \to \psi \To \Delta_1, \Delta_2)$ is obtained by
            the \SequentRuleName{\to}{L} rule,
            then by the induction hypothesis, both $\Gamma_1^\flat \To \varphi^\sharp, \Delta_1^\sharp$ and $\Gamma_2^\flat, \psi^\flat \To \Delta_2^\sharp$ are provable in $\LK$.
            Here, $(\varphi \to \psi)^\flat = \varphi^\sharp \to \psi^\flat$. Then:
            \begin{center}
              \Axiom$\Gamma_1^\flat \fCenter \varphi^\sharp, \Delta_1^\sharp$
              \Axiom$\Gamma_2^\flat, \psi^\flat \fCenter \Delta_2^\sharp$
              \RightLabel{\scriptsize{(\SequentRuleName{\to}{L})}}
              \BinaryInf$\Gamma_1^\flat, \Gamma_2^\flat, \varphi^\sharp \to \psi^\flat \fCenter \Delta_1^\sharp, \Delta_2^\sharp$
              \DisplayProof
            \end{center}
      \item If $(\Gamma \To \Delta) = (\Gamma \To \varphi \to \psi, \Delta')$ is obtained by
            the \SequentRuleName{\to}{R} rule,
            then $\Gamma^\flat, \varphi^\flat \To \psi^\sharp, \Delta'^\sharp$ is provable by the induction hypothesis.
            Here, $(\varphi \to \psi)^\sharp = \varphi^\flat \to \psi^\sharp$. Then:
            \begin{center}
              \Axiom$\Gamma^\flat, \varphi^\flat \fCenter \psi^\sharp, \Delta'^\sharp$
              \RightLabel{\scriptsize{(\SequentRuleName{\to}{R})}}
              \UnaryInf$\Gamma^\flat \fCenter \varphi^\flat \to \psi^\sharp, \Delta'^\sharp$
              \DisplayProof
            \end{center}
    \end{itemize}

    Now we consider the cases where $\Gamma \To \Delta$ is obtained by one of the new rules in $\GNxAmn$.

    If $(\Gamma \To \Delta) = (\To \Box \varphi)$ is obtained by $\frac{\To \varphi}{\To \Box \varphi}$ (nec),
    then $\NxAmn \vdash \varphi$, which implies that $(\Box \varphi)^\sharp = \top$, and $\To \top$ is provable in $\LK$.

    If $\NxAmn = \NpAmn$, $m = 0$, $n \ge 2$, and $(\Gamma \To \Delta) = (\Box \Box \varphi \To)$ is obtained by $\frac{\Box \varphi \To}{\Box \Box \varphi \To}$ (rosbox),
    then $\NpAmn \vdash \neg \Box \varphi$, which implies that $(\Box \Box \varphi)^\flat = \bot$, and $\bot \To$ is an initial sequent in $\LK$.

    If $\NxAmn = \NRAmn$ and $(\Gamma \To \Delta) = (\Box \varphi \To)$ is obtained by $\frac{\varphi \To}{\Box \varphi \To}$ (ros),
    then $\NRAmn \vdash \neg \varphi$, which implies that $(\Box \varphi)^\flat = \bot$, and $\bot \To$ is an initial sequent in $\LK$.

    If $n > m$ and $(\Gamma \To \Delta) = (\Box^n \varphi, \Gamma' \To \Delta)$ is obtained by $\frac{\Box^m \varphi, \Box^n \varphi, \Gamma' \To \Delta}{\Box^n \varphi, \Gamma' \To \Delta}$ (accL),
    then by the induction hypothesis, $(\Box^m \varphi)^\flat, (\Box^n \varphi)^\flat, \Gamma'^\flat \To \Delta^\sharp$ is provable in $\LK$.
    Here, there are $k \ge n$ and $\psi \in \MF$ such that $\Box^n \varphi = \Box^k \psi$ and $\psi \ne \Box \chi$ for any $\chi \in \MF$.
    We distinguish the following cases.
    \begin{itemize}
      \item If $(\Box^n \varphi)^\flat = (\Box^k \psi)^\flat = \bot$, then:
            \begin{center}
              \Axiom$\bot \fCenter $
              \RightLabel{\scriptsize{(wL, wR)}}
              \UnaryInf$\bot, \Gamma'^\flat \fCenter \Delta^\sharp$
              \DisplayProof
            \end{center}
      \item Otherwise, $k \ge n > m$ implies $(\Box^n \varphi)^\flat = (\Box^k \psi)^\flat = p_{\Box^k \psi} \land (\Box^{k-n+m} \psi)^\flat = p_{\Box^n \varphi} \land (\Box^m \varphi)^\flat$. Then:
            \begin{center}
              \scalebox{0.8}{
                \Axiom$(\Box^m \varphi)^\flat \fCenter (\Box^m \varphi)^\flat$
                \RightLabel{\scriptsize{(\SequentRuleName{\land}{L})}}
                \UnaryInf$p_{\Box^n \varphi} \land (\Box^m \varphi)^\flat \fCenter (\Box^m \varphi)^\flat$

                \Axiom$(\Box^m \varphi)^\flat, p_{\Box^n \varphi} \land (\Box^m \varphi)^\flat, \Gamma'^\flat \fCenter \Delta^\sharp$

                \RightLabel{\scriptsize{(cut)}}
                \BinaryInf$p_{\Box^n \varphi} \land (\Box^m \varphi)^\flat, \Gamma'^\flat \fCenter \Delta^\sharp$
                \DisplayProof
              }
            \end{center}
    \end{itemize}

    If $m > n$ and $(\Gamma \To \Delta) = (\Gamma \To \Delta', \Box^m \varphi)$ is obtained by $\frac{\Gamma \To \Delta', \Box^m \varphi, \Box^n \varphi}{\Gamma \To \Delta', \Box^m \varphi}$ (accR),
    then by the induction hypothesis, $\Gamma^\flat \To \Delta'^\sharp, (\Box^m \varphi)^\sharp, (\Box^n \varphi)^\sharp$ is provable in $\LK$.
    Here, there are $k \ge m$ and $\psi \in \MF$ such that $\Box^m \varphi = \Box^k \psi$ and $\psi \ne \Box \chi$ for any $\chi \in \MF$.
    We distinguish the following cases.
    \begin{itemize}
      \item If $(\Box^m \varphi)^\sharp = (\Box^k \psi)^\sharp = \top$, then:
            \begin{center}
              \Axiom$\bot \fCenter \bot$
              \RightLabel{\scriptsize{(\SequentRuleName{\to}{R})}}
              \UnaryInf$ \fCenter \top$
              \RightLabel{\scriptsize{(wL, wR)}}
              \UnaryInf$\Gamma^\flat \fCenter \Delta'^\sharp, \top$
              \DisplayProof
            \end{center}
      \item Otherwise, $k \ge m > n$ implies $(\Box^m \varphi)^\sharp =  (\Box^k \psi)^\sharp = p_{\Box^k \psi} \lor (\Box^{k-m+n} \psi)^\sharp = p_{\Box^m \varphi} \lor (\Box^n \varphi)^\sharp$. Then:
            \begin{center}
              \scalebox{0.8}{
                \Axiom$\Gamma^\flat \fCenter \Delta'^\sharp, p_{\Box^m \varphi} \lor (\Box^n \varphi)^\sharp, (\Box^n \varphi)^\sharp$

                \Axiom$(\Box^n \varphi)^\sharp \fCenter (\Box^n \varphi)^\sharp$
                \RightLabel{\scriptsize{(\SequentRuleName{\lor}{R})}}
                \UnaryInf$(\Box^n \varphi)^\sharp \fCenter p_{\Box^m \varphi} \lor (\Box^n \varphi)^\sharp$

                \RightLabel{\scriptsize{(cut)}}
                \BinaryInf$\Gamma^\flat \fCenter \Delta'^\sharp, p_{\Box^m \varphi} \lor (\Box^n \varphi)^\sharp$
                \DisplayProof
              }
            \end{center}
    \end{itemize}
    In any case, $(\Gamma^\flat \To \Delta^\sharp)$ is provable in $\LK$.
  \end{proof}
\end{theorem}

\begin{corollary} [2.\ of Def.\ \ref{def:p18n}] \label{cor:embedding}
  $\NxAmn \vdash \varphi \to \psi$ implies $\CL' \vdash \varphi^\flat \to \psi^\sharp$.
\end{corollary}

\begin{remark}
  We shall note that having a cut-free proof of $\Gamma \To \Delta$ is crucial in proving Theorem \ref{thm:emulation},
  although it does not seem to rely on the cut elimination of $\GNxAmn$ at first glance.
  Suppose that $(\Gamma \To \Delta) = (\Gamma_1, \Gamma_2 \To \Delta_1, \Delta_2)$ were obtained by the cut rule
  $\left(\frac{\Gamma_1 \To \Delta_1, \varphi \quad \varphi, \Gamma_2 \To \Delta_2}{\Gamma_1, \Gamma_2 \To \Delta_1, \Delta_2}\right)$,
  then by the induction hypothesis, $\Gamma_1^\flat \To \Delta_1^\sharp, \varphi^\sharp$ and $\varphi^\flat, \Gamma_2^\flat \To \Delta_2^\sharp$ would be provable in $\LK$.
  As $\varphi^\sharp$ is provably weaker than $\varphi^\flat$, there would be no way of applying the cut rule to these two sequents and thus obtaining $\Gamma_1^\flat, \Gamma_2^\flat \To \Delta_1^\sharp, \Delta_2^\sharp$.
\end{remark}

\begin{theorem}
  $(\sharp, \flat)$ is a propositionalization of $\NxAmn$ into $\CL$.
  \begin{proof}
    Follows from Propositions \ref{prop:sharpflat-variables}, \ref{prop:equality} and Corollary \ref{cor:embedding}.
  \end{proof}
\end{theorem}

\begin{corollary}
  $\NxAmn$ enjoys ULIP.
  \begin{proof}
    Follows from Theorem \ref{thm:p18n-ulip}.
  \end{proof}
\end{corollary}

\section*{Acknowledgements}
\addcontentsline{toc}{section}{Acknowledgements}

The author would like to express our deep appreciation for Taishi Kurahashi, the supervisor, who provided us the extensive support and advice throughout this research.
The author would also like to thank Professor Katsuhiko Sano for giving us some helpful comments in the earlier stages of this research,
and the anonymous referees for giving valuable advice on improving this paper.

\bibliographystyle{unsrt}
\bibliography{refs}

\end{document}